\documentclass{amsart}

\usepackage{color}
\usepackage{graphicx,epsf}
\newtheorem{theorem}{Theorem}[section]
\newtheorem{lemma}[theorem]{Lemma}

\newtheorem{proposition}[theorem]{Proposition}
\theoremstyle{definition}

\theoremstyle{remark}
\newtheorem{remark}[theorem]{Remark}

\numberwithin{equation}{section}

\newcommand{\abs}[1]{\lvert#1\rvert}

\newcommand{\blankbox}[2]{%
  \parbox{\columnwidth}{\centering
    \setlength{\fboxsep}{0pt}%
    \fbox{\raisebox{0pt}[#2]{\hspace{#1}}}%
  }%
}
\newcommand {\sign} {\mbox{\bf sign}}

\newcommand {\OO} {{\mathcal O}}
\newcommand {\R} {{\rm R}}
\newcommand {\C} {{\rm C}}

\newcommand {\eop}      {\hfill $\Box$}
\newcommand {\junk}[1]{}
\newcommand {\hide}[1]{}
\newcommand {\comment}[1]{\textcolor{red}{{\sc Comment: #1}}}

\newcommand {\s}        {\mbox{\rm sign}}
\newcommand {\sz}       {\mbox{\rm zero}}
\newcommand {\D}     {\mbox{\rm D}}
\newcommand {\K}     {\mbox{\rm K}}
\newcommand {\A}     {\mbox{\rm A}}
\newcommand{\LL}     {\mbox{\rm L}}
\newcommand {\cf}     {\mbox{\rm cof}}
\newcommand {\SF}     {\mbox{\rm SF}}
\newcommand {\rem}     {\mbox{\rm Rem}}
\newcommand {\quo}     {\mbox{\rm Quo}}
\newcommand {\prem}     {\mbox{\rm Prem}}
\newcommand {\pquo}     {\mbox{\rm Pquo}}
\newcommand {\pagcd}    {\mbox{\rm Posgcd}}
\newcommand {\paremseq} {\rm TRems}
\newcommand {\parem}   {\mbox{\rm Parem}}
\newcommand {\posquo}   {\mbox{\rm Posquo}}
\newcommand {\posgcd}    {\mbox{\rm Posgcd}}
\newcommand {\lc}       {\mbox{\rm lcof}}
\newcommand {\Real}[1]   {\mbox{${\mathbb R}^{#1}$}}
\newcommand {\Complex}[1]   {\mbox{${\mathbb C}^{#1}$}}
\newcommand {\Rd}        {\Real{d}}
\newcommand {\Rdn}       {\Real{dn}}
 \newcommand {\Rtwo}      {\Real{2}}
\newcommand {\Sphere}[1] {\mbox{${\bf S}^{#1}$}}
\newcommand {\Rm}        {\Real{m}}
 \newcommand {\re}         {\Real{}}
 \newcommand {\complex}   {\Complex{}}
\newcommand {\SI}         {\mbox{\rm SIGN}_R}
\newcommand {\Sign}         {\mbox{\rm SIGN}}
\newcommand {\level}{\mbox{\rm level}}
\newcommand {\Z}  {{\mathbb Z}}
 \newcommand {\N}         {{\mathbb N}}
\newcommand {\Q}         {{\mathbb Q}}

\newcommand{\Qu} {\rm Q}
\newcommand{\F}{\mbox{\rm F}}
\newcommand{\decide} {\rm Decide}

\newcommand  {\pos} {{\rm pos}}
\newcommand {\ZZ} {{\rm Z}}
\newcommand {\RR} {{\mathcal R}}
\newcommand {\V} {{\rm V}}
\newcommand {\Der} {{\rm Der}}
\newcommand {\SQ} {{\rm SQ}}

\newcommand {\la}   {{\langle}}
\newcommand {\ra}   {{\rangle}}
\newcommand {\eps} {{\varepsilon}}
\newcommand {\E} {{\rm Ext}}
\newcommand {\dist} {{\rm dist}}

\newcommand {\Syl} {\mbox{\rm Syl}}
\newcommand {\pdet} {\mbox{\rm pdet}}
\newcommand {\I} {\mbox{\rm I}}
\newcommand {\Id} {\mbox{\rm Id}}
\newcommand {\SH} {\mbox{\mathcal SH}}
\newcommand {\SR} {\mbox{\rm SR}}
\newcommand {\sr} {\mbox{\rm sr}}
\def\gau{{\rm Gau}}
\def\sign{{\rm sign}}
\def\han{{\rm Han}}
\newcommand {\lcm}      {\mbox{\rm lcm}}
\newcommand {\ldp}      {\mbox{\rm ldp}}
\newcommand {\GRem}      {\mbox{\rm GRem}}
\newcommand {\Res}      {\mbox{\rm Res}}

\newcommand {\Ker}      {\mbox{\rm Ker}}
\newcommand {\Ima}      {\mbox{\rm Im}}

\newcommand {\lmon}      {\mbox{\rm lmon}}
\newcommand{\Sy}  {\mbox{\rm S}}
\newcommand {\lt}      {\mbox{\rm lt}}
\newcommand {\PP}     {{\mathbb P}}

\newcommand {\SU} {\mbox{\rm SU}}
\newcommand {\SV} {\mbox{\rm SV}}
\newcommand {\SB} {\mbox{\rm SB}}

\newcommand {\W}      {\mbox{\rm W}}
\newcommand {\T}      {\mbox{\rm T}}
\newcommand {\Rec}      {{\rm Rec}}
\newcommand {\Ho}      {\mbox{\rm H}}
\newcommand {\Co}      {\mbox{\rm C}}
\newcommand {\M}        {\mbox{\rm M}}
\newcommand {\Le}       {\mbox{\rm L}}

\newcommand{\EQ}  {\mbox{\rm EQ}}
\newcommand{\epc}  {\mbox{\rm EPC}}

\newcommand {\NF}       {\mbox{\rm NF}}
\newcommand {\MUS}      {\mbox{\rm U}}
\newcommand {\B} {\mbox {\rm B}}
\newcommand {\Red} {\mbox {\rm Red}}

\def\addots{\mathinner{\mkern1mu
\raise1pt\vbox{\kern7pt\hbox{.}}
\mkern2mu\raise4pt\hbox{.}\mkern2mu
\raise7pt\hbox{.}\mkern1mu}}

\newcommand {\Def} {\mbox {\rm Def}}
\newcommand {\Cr} {\mbox {\rm Cr}}

\newcommand{\Po}  {\mbox{\rm Pol}}
\newcommand{\RM}  {\mbox{\rm RM}}
\newcommand{\URM}  {\mbox{\rm URM}}

\newcommand{\pr}  {\mbox{\rm Proj}}
\renewcommand {\Im} {\mbox {\rm Im}}

\newcommand{\HH}  {{\mathrm H}}

\newcommand{\coucou}[1]{\ifvmode\else\marginpar[\hfill$\rhd$]{$\lhd$}\fi
                       $\langle$\textsc{#1}$\rangle$}
\newcommand{\coucoux}[1]{$\langle$\textsc{#1}$\rangle$}
\newcommand{\mf}[1]{\coucou{\underline{MF}: #1}}
\newcommand{\w}[1]{\coucou{\underline{WARNING}: #1}}

\newcommand{\basu}[1]{}
\newcommand {\rp}[1]{}

\begin{document}

\title
[
A tight bound on the number
of connected components
]
{
An asymptotically tight bound on the number of semi-algebraically connected components of 
realizable sign conditions
\footnote{2000 Mathematics Subject Classification 14P10, 14P25}}

\author{Saugata Basu}
\address{School of Mathematics,
Georgia Institute of Technology, Atlanta, GA 30332, U.S.A.}
\email{saugata.basu@math.gatech.edu}
\thanks{Supported in part by an NSF Career Award 0133597 and
an Alfred P. Sloan Foundation Fellowship.}

\author{Richard Pollack}
\address{Courant Institute of Mathematical Sciences, 
New York University, New York, NY 10012, U.S.A.}
\email{pollack@cims.nyu.edu}
\thanks{Supported in part by 
NSA  grant MDA904-01-1-0057 and NSF grants
CCR-9732101 and  CCR-0098246.}

\author{Marie-Fran\c{c}oise Roy}
\address{IRMAR (URA CNRS 305), 
Universit\'{e} de Rennes,
Campus de Beaulieu 35042 Rennes cedex FRANCE.}
\email{marie-francoise.roy@univ-rennes1.fr}
%    General info
%%\subjclass{Primary 14P10, 14P25; Secondary 46E25, 20C20}

%%\date{January 1, 1994 and, in revised form, June 22, 1994.}

%%\dedicatory{This paper is dedicated to our authors.}


\begin{abstract}
We prove an asymptotically  tight bound 
(asymptotic with respect to the number of polynomials 
for fixed degrees and number of variables) 
on the 
number of semi-algebraically connected components of the realizations of all realizable sign
conditions of a family of real polynomials.
More precisely, we prove that the number of semi-algebraically connected components of the
realizations of all realizable sign conditions of a family of $s$ polynomials
in $\R[X_1,\ldots,X_k]$ whose degrees are at most $d$ is bounded by
\[
\frac{(2d)^k}{k!}s^k + O(s^{k-1}).
\]
This improves the best upper bound known previously which was
\[
\frac{1}{2}\frac{(8d)^k}{k!}s^k + O(s^{k-1}).
\]
The new bound matches asymptotically the lower
bound obtained for families of polynomials each of which is a product
of generic polynomials of degree one.
\end{abstract}

\maketitle

\section{Introduction}
Let $\K$ be a field and ${\mathcal P} \subset \K[X_1,\ldots,X_k]$ be a 
finite family of polynomials with $\#{\mathcal P} = s$. 
Then ${\mathcal P}$ naturally induces a 
partition of $\K^k$ into constructible subsets, such that over 
each subset belonging to the partition, each polynomial 
$P \in {\mathcal P}$  either vanishes at every point of the subset or 
is non-zero at every point of the subset 
and each such set is maximal with respect to this property.
If moreover  $\K$ is an ordered field, we can consider a 
refined partition, such that each $P \in {\mathcal P}$ maintains its 
{\em sign} over each element of the partition. 
Notice that the number of sets in the partition can be at most $2^{s}$ in 
the unordered and $3^{s}$ in the ordered case. However, 
much tighter bounds are known in both cases 
\cite{BPR95a,RBG01}. 
In fact, in case the field $\K$ is ${\mathbb R}$ or ${\mathbb C}$, 
the number of semi-algebraically connected components
of the partition of $\K^k$ induced by ${\mathcal P}$ is bounded by a polynomial
function of $s$ of degree $k$. 

In this paper we prove asymptotically tight bounds on the number of 
semi-algebraically connected components of sets in  
the partitions described above. Our bounds are asymptotic in 
$s$, with the number of variables and the degrees of
the polynomials in ${\mathcal P}$ considered fixed. 
Asymptotics with respect to the number of polynomials (with the degrees and
the number of variables considered fixed) is considered important in 
applications in discrete and computational geometry (see \cite{Matousek}).

The proofs of all previous results 
stating that the number 
of semi-algebraically connected components
of the partition of $\K^k$ induced by ${\mathcal P}$ is bounded by a polynomial
function of $s$ of degree $k$
(see \cite{Warren, Alon, BPR95a,BPR8}),
used a common technique of replacing the given family of polynomials,
${\mathcal P}$, 
by another family, ${\mathcal P}^{\star}$,
obtained by infinitesimal perturbations of polynomials in
${\mathcal P}$. 
The partition induced by the family 
${\mathcal P}^{\star}$ is closely related to that
of ${\mathcal P}$ and at the same time 
the family ${\mathcal P}^{\star}$ has useful properties 
which makes it easier to bound the topological
complexity of the partition induced by it.
However, the cardinality of ${\mathcal P}^{\star}$ is larger (often 
by a factor of $4$)
than that of ${\mathcal P}$, 
and this fact introduces an extra factor of $4^k$ in the
upper bound. 
In this paper we use a new technique 
(see Section \ref{sec:summary} below) that avoids replacing 
${\mathcal P}$ by a larger family of polynomials and thus
we are able to obtain a tight asymptotic bound.

\section{Preliminaries}
\label{sec:preliminaries}
We begin  with a few definitions.

Let $\R$ be a real closed field and $\C$ its algebraic closure. 
For $x \in \R$ we define
\[\s(x) =
\begin{cases}
0& {\rm if\ 
}\ x=0,\cr
1& {\rm if\ 
}\ x> 0,\cr
-1& {\rm if\ 
}\ x < 0.
\end{cases}
\]
Similarly for $x \in \C$ we define
\[\sz(x) =
\begin{cases}
0& {\rm if\ 
}\ x=0,\cr
1& {\rm if\ 
}\ x\neq 0.
\end{cases}
\]

Let ${\mathcal P}$ be a finite subset of $\R[X_1,\ldots,X_k]$.
A  {\em sign condition}  (resp. {\em zero pattern}) on
${\mathcal P}$ is an element of $\{0,1,- 1\}^{\mathcal P}$
(resp. $\{0,1\}^{\mathcal P}$). 
Given ${\mathcal P}$ we fix $\Omega > 0$ to be a sufficiently large
element of $\R$ and denote by

$$
\displaylines{
B_k(0,\Omega) = \{ x \in \R^k \mid |x|^2 < \Omega\}, \\
\overline{B_k(0,\Omega)} = \{ x \in \R^k \mid |x|^2 \leq \Omega\}.
}
$$

In our arguments we will often consider subsets
of $\C^k \cong \R^{2k}$ and we will denote the open ball centered at the
origin of radius $\Omega$ in $\C^k = \R^{2k}$ 
also by $B_{2k}(0,\Omega)$,  that is 
$$
\displaylines{
B_{2k}(0,\Omega) = \{ x \in \C^k \mid |x|^2 < \Omega \}, \\
\overline{B_{2k}(0,\Omega)} = \{ x \in \C^k \mid |x|^2 \leq \Omega \}.
}
$$

We will restrict our attention to realizations of sign conditions (resp.
zero patterns) inside the 
closed
ball $\overline{B_k(0,\Omega)}$ (resp. $\overline{B_{2k}(0,\Omega)}$).

The {\em realization} of the sign condition $\sigma$
is the semi-algebraic set

\begin{equation}
\label{eqn:R(Z)}
\RR(\sigma,\R^k)= \{x\in \R^k\;\mid\; 
\bigwedge_{P\in{\mathcal P}} \s({P}(x))=\sigma(P) \} 
\cap \overline{B_{k}(0,\Omega)}.
\end{equation}

We call a semi-algebraically connected component $C$ of $\RR(\sigma,\R^k)$ 
bounded if
$C \cap \partial \overline{B_{k}(0,\Omega)} = \emptyset$ and unbounded otherwise.

The {\em realization} of the   zero pattern $\rho$
is the  set 
\begin{equation}
\label{eqn:C(Z)}
\RR(\rho,\C^k)   =      \{x\in \C^k\;\mid\; 
\bigwedge_{P\in{\mathcal P}} \sz({P}(x))=\rho(P) \}
\cap \overline{B_{2k}(0,\Omega)}.
\end{equation}

It is a consequence of Hardt's triviality theorem 
\cite{Hardt}
that the homeomorphism
type of $\RR(\sigma,\R^k)$ (resp. $\RR(\rho,\C^k)$) is constant for all sufficiently
large $\Omega > 0$.

We 
denote
the set of zeros
of ${\mathcal P}$ in $\R^k$  (resp. in $\C^k$) by
$$
\displaylines{
\ZZ({\mathcal P},\R^k)=\{x\in \R^k\mid\bigwedge_{P\in{\mathcal P}}P(x)= 0\}
}
$$

$$
\displaylines{
(\mbox{resp. }
\ZZ({\mathcal P},\C^k)=\{x\in \C^k\mid\bigwedge_{P\in{\mathcal P}}P(x)= 0\}
).
}
$$ 

We denote by
$$
\displaylines{
{\rm Sign}({\mathcal P}) = 
\{\sigma \in \{0,1,-1\}^{\mathcal P}\mid \RR(\sigma,\R^k) \neq \emptyset\},
}
$$
and 
$$
\displaylines{
\mbox{Zero-pattern}({\mathcal P})=
\{\rho \in \{0,1\}^{\mathcal P}\mid \RR(\rho,\C^k) \neq \emptyset\}.
}
$$

For $\sigma\in \{0,1,-1\}^{\mathcal P}$ 
we will denote by
$b_i^{BM}(\sigma,\R^k)$ the dimension of 

\[
\HH_i^{BM}(\RR(\sigma,\R^k)),
\] 
the $i$-th Borel-Moore homology group of the locally closed set
$\RR(\sigma,\R^k)$ with coefficients in $\Z_2$.
Similarly, we will denote by
$b_i(\sigma,\R^k)$ the dimension of 
$\HH_i(\RR(\sigma,\R^k))$, 
the $i$-th singular homology group of 
$\RR(\sigma,\R^k)$ with coefficients in $\Z_2$.
We will denote by
\[
b^{BM}(\sigma,\R^k) = \sum_{i \geq 0} b_i^{BM}(\sigma,\R^k).
\]

Similarly,
for $\rho \in \{0,1\}^{\mathcal P}$ 
we will denote by
$b_i^{BM}(\rho,\C^k)$ the dimension of 
\[
\HH_i^{BM}(\RR(\rho,\C^k)),
\] 
the $i$-th Borel-Moore homology group of the locally closed set
$\RR(\rho,\C^k)$ with coefficients in $\Z_2$.
We will denote by
\[
b^{BM}(\rho,\C^k) = \sum_{i \geq 0} b_i^{BM}(\rho,\C^k).
\]

More generally,
for any locally closed semi-algebraic set $X$ we denote by
$b_i^{BM}(X)$ to be 
the dimension of $\HH_i^{BM}(X)$ with $\Z_2$-coefficients,
and we will denote by
\[
b^{BM}(X)=  \sum_{i \geq 0} b_i^{BM}(X).
\]

It will be also useful to have following notation: 
$$
\displaylines{
b_{i,\R}(s,d,k) = \max_{{\mathcal P}} b_i({\mathcal P},\R^k),
}
$$
where
$$
\displaylines{
b_i({\mathcal P},\R^k) = \sum_{\sigma \in \{0,1,-1\}^{\mathcal P}}b_i(\sigma,\R^k),
}
$$
and the maximum is taken over all
${\mathcal P} \subset \R[X_1,\ldots,X_k]$
with $ \#{\mathcal P} = s$ and  $\deg(P) \leq d$ for each  $P \in {\mathcal P}$.

In some applications, 
for instance in bounding the number of combinatorial types of polytopes
or order types of point configurations (see \cite{GP2,GP86}),
one is 
only interested in the number of realizable sign conditions
or zero patterns and not in any topological properties of their
realizations. In these situations it is useful to obtain bounds on
$$
\displaylines{
N_{\R}(s,d,k) = \max_{{\mathcal P}} \#{\rm Sign}({\mathcal P}), 
}
$$
where the maximum is taken over all
${\mathcal P} \subset \R[X_1,\ldots,X_k]$
with $ \#{\mathcal P} = s$ and  $\deg(P) \leq d$ for each  $P \in {\mathcal P}$,
as well as
$$
\displaylines{
N_{\C}(s,d,k) = \max_{{\mathcal P}} \#{\mbox{Zero-pattern}}({\mathcal P}),
}
$$
where the maximum is taken over all
${\mathcal P} \subset \C[X_1,\ldots,X_k]$
with $ \#{\mathcal P} = s$ and  $\deg(P) \leq d$ 
for each  $P \in {\mathcal P}$.
Note that 
$N_{\R}(s,d,k)$ is the maximum number of sign conditions realized in $\R^k$ by
a set of polynomials of size $s$ and degrees bounded by $d$.

\section{Known Bounds}
\subsection{Zero patterns}
The problem of bounding the number of realizable zero patterns 
of a sequence of  $s$
polynomials in $k$ variables of degree at most $d$ 
over an arbitrary field was considered by Ronyai et al. in \cite{RBG01},
where they prove 
\begin{equation}
\label{eqn:rbg}
N_{\C}(s,d,k) \leq  {sd+k \choose k} 
\leq \left(\frac{d^k}{k!}\right) s^k + O(s^{k-1}) 
\end{equation}
using linear algebra arguments. 
However, their method is not useful for the 
problem of bounding the number of realizable sign conditions over an ordered
field.
%%sb adds
A slightly better bound (see Theorem 5.2 below)
obtained as a consequence of the methods used in this paper 
was proved earlier by Jeronimo and Sabia in \cite{JS00}.

\subsection{Sign conditions}
We now consider the case of sign conditions.
When $d=1$, we have 
$N_{\R}(s,1,k) = b_{0,\R}(s,1,k)$, since the realization of each realizable
sign condition in this case is a convex polyhedron, which is clearly semi-algebraically connected.
It is also easy in this case to deduce an exact expression for 
$N_{\R}(s,1,k) = b_{0,\R}(s,1,k)$,
namely,
\begin{eqnarray}
\label{eqn:linear}
N_{\R}(s,1,k) &=& b_{0,\R}(s,1,k) \nonumber\\
&=& 
\sum_{i=0}^{k} \sum_{j=0}^{k-i}{s \choose i}{s-i \choose j}\nonumber\\
&=&
\left(\sum_{i=0}^{k} \frac{1}{i!(k-i)!}\right)s^k + O(s^{k-1}) \nonumber \\
&=& \frac{2^k}{k!} s^k + O(s^{k-1}).
\end{eqnarray}

This is the number of cells of all dimensions in an arrangement of $s$
hyperplanes in general position in $\R^k$
and hence this bound is sharp. 
For $0 \ll k \ll s$ we have
\[N_{\R}(s,1,k)= b_{0,\R}(s,1,k) \sim \left(\frac{2es}{k}\right)^k
\]
using Stirling's approximation for $k!$.

For $d >1$ Alon \cite{Alon} proved 
a bound of 
\[
\left(\frac{8esd}{k}\right)^k
\]
on $N_{\R}(s,d,k)$. 
Previously, Warren \cite{Warren} had proved a bound of 
\[
\left(\frac{4esd}{k}\right)^k
\]
on the number of realizable strict sign conditions (that is sign conditions
$\sigma$ such that $\sigma(P) \neq 0$ for all $P \in {\mathcal P}$).  

The following bounds on the sums of the
individual Betti numbers of the realizations of all realizable sign conditions
restricted to a real variety of real dimension $k'$ also defined by 
polynomials of degree at most $d$ was proved in \cite{BPR8}.

\begin{eqnarray}
\label{eqn:pams}
        b_{i,\R}(s,d,k,k') &\le &\sum_{j=0}^{k' - i} 
              {s \choose j} 4^{j}  d(2d-1)^{k-1} \nonumber\\
&= &\left(\frac{2^{2k'+k-2i -1}d^k}{(k'-i)!}\right)s^{k'-i} + O(s^{k'-i-1})
\end{eqnarray}
(where the last parameter $k'$ denotes the dimension of the ambient 
variety).

Note that $b_{0,\R}({\mathcal P})$
is the total number of semi-algebraically connected components
of the realizations of all  realizable sign conditions of
${\mathcal P}$ and the above result gives an asymptotic
bound of 
\[
\frac{1}{2}\left(\frac{8esd}{k}\right)^k
\]
(which is asymptotically same as that proved by Alon in \cite{Alon})
on the number of semi-algebraically connected components of the realizations of all
realizable sign conditions which is a priori larger than just the
number of realizable sign conditions. This distinction is 
important since in many applications,
for instance in bounding the number of isotopy classes  of point
configurations,
it is the number of semi-algebraically connected components which is 
important.

In this paper we consider the problem of proving an asymptotic
upper bound on $b_{0,\R}(s,d,k)$ (that is the number of semi-algebraically connected components
of all realizable sign conditions of a family of $s$ polynomials
in $\R[X_1,\ldots,X_k]$ of degrees bounded by $d$), 
for fixed $d,k$ and large $s$.
It follows from (\ref{eqn:pams}) that
$b_{0,\R}(s,d,k)$ is bounded from above by a polynomial in $s$ of degree $k$, 
namely,
$$
\displaylines{
        b_{0,\R}(s,d,k) \le 
\left(\frac{2^{3k - 1}d^k}{k!}\right)s^{k} + O(s^{k-1}).
}
$$
The leading coefficient of the above bound is 
$
\sim \left(\frac{8ed}{k}\right)^k.
$

On the other hand by taking $s$ polynomials each a product of $d$
generic linear polynomials, we have by (\ref{eqn:linear}) a lower bound of 
\begin{eqnarray*}
b_{0,\R}(s,d,k) &\geq& \sum_{i=0}^{k} \sum_{j=0}^{k-i}{sd \choose i}{sd-i \choose j}\\
&=& 
\left(\sum_{i=0}^k \frac{d^k}{i! (k -i)!} \right)s^k +  \Theta(s^{k-1})\\
&=& \frac{(2d)^k}{k!} s^k +  \Theta(s^{k-1}).
\end{eqnarray*}

The leading coefficient of the lower bound is,
\[
\frac{(2d)^k}{k!} \sim \left(\frac{2ed}{k}\right)^k.
\]
(In the above bounds we provide the asymptotic estimates on the right for
aid in comparison.)

Clearly the upper and lower estimates differ by a factor of $4^k$.
As mentioned in the introduction,
this gap is due to the fact that the technique used in proving the 
upper bound (\ref{eqn:pams}) involves replacing the given family of polynomials
by another family obtained by infinitesimal perturbations of the original
polynomials. The realizations of a certain subset of {\em strict} 
sign conditions 
of this new family are then shown to be in one to one correspondence with
the realizations of all realizable sign conditions of the original family.
Moreover, the corresponding realizations have the same homotopy type. 
However, the cardinality of the new family is four times that of the
original family, and this fact introduces an extra factor of $4^k$ in the
upper bound.

\section{Brief Summary of our method}
\label{sec:summary}
In this paper we use a new technique for proving an upper
bound on $b_{0,\R}(s,d,k)$ that avoids replacing the given
family of polynomials by a larger family of polynomials.
The new idea is to consider the {\em zero patterns} of the given family of 
polynomials over $\C^k$, where we are able to use degree theory
of complex varieties and certain generalized Bezout inequalities,
as well as inequalities derived from the Mayer-Vietoris exact sequence, 
to directly obtain an asymptotically tight bound on the sum of the 
Borel-Moore Betti numbers of the realizations of all realizable 
zero patterns of the original family of polynomials in $\C^k$. 
This in turn provides an asymptotically tight bound on the sum of the 
Borel-Moore Betti numbers of the  realizations of all 
realizable sign conditions
of the same family in $\R^k$ using Smith inequalities
(see Proposition \ref{prop:smith} below). Note that even though we are
interested in bounding the number of semi-algebraically connected components of sign 
conditions over  $\R^k$,
in order to use Smith inequality we need to bound the higher Betti numbers
of the zero patterns in $\C^k$. 

The reason for using Borel-Moore homology instead of singular homology is
that Borel-Moore homology has an useful additivity property 
(see Proposition \ref{prop:BM} below)
not satisfied by singular homology. This property plays an important role
in our proofs. 
However, in order to obtain bounds on the number of semi-algebraically connected components of
the realizations of sign conditions,
we need to relate the Borel-Moore Betti numbers
of the realizations of sign conditions, with the number of semi-algebraically connected
components. We prove that the number of bounded semi-algebraically connected components of any
fixed sign conditions is bounded by the sum of the Borel-Moore Betti numbers
of the realization of the sign condition
(see Lemma \ref{lem:2} below). This allows us to bound the
number of bounded semi-algebraically connected components in terms of the Borel-Moore Betti
numbers. We bound the number of unbounded components separately using
Inequality \eqref{eqn:pams}.

\section{Main Results}
\label{sec:results}
We obtain the following bound on the number of semi-algebraically connected components
of the realizations of all realizable sign conditions of a family of 
polynomial.

\begin{theorem}
\label{the:main1}
\begin{eqnarray*}
b_{0,\R}(s,d,k) &\leq& 
\sum_{0 \leq \ell \leq k}
\left({s \choose k - \ell}{s  \choose \ell}d^k + 
\sum_{1 \leq i \leq \ell} {s \choose k - \ell}{s \choose \ell - i}d^{O(k^2)}\right)\\
&=& \sum_{\ell=0}^k {s \choose k - \ell}{s \choose \ell}d^k + O(s^{k-1})\cr
&=& \left(\sum_{\ell=0}^k \frac{d^k}{\ell! (k -\ell)!} \right)s^k + O(s^{k-1}) \cr
&=& \frac{(2d)^k}{k!} s^k + O(s^{k-1}).
\end{eqnarray*}
\end{theorem}

For $0 < d , k \ll s$, this gives
\[
b_{0,\R}(s,d,k) \sim 
\left(\frac{2esd}{k}\right)^k.
\]
This matches  asymptotically 
the lower bound obtained by taking $s$ polynomials each of which
is a product of $d$ linear polynomials in general position.

In the process of proving Theorem \ref{the:main1}, we slightly improve 
the bound on $N_{\C}(s,d,k)$ proved by Ronyai et al.  
\cite{RBG01} mentioned above. 
%%sb adds
This result was already 
obtained previously in \cite{JS00}.

Denoting by $N_{\C,_\ell}(s,d,k)$ the maximum number of realizable
zero patterns whose realizations are of (complex) dimension $\ell$, we have
\begin{theorem}
\label{the:0-1}
\[
N_{\C,\ell}(s,d,k) \leq {s \choose k-\ell} d^{k - \ell}.
\]
\end{theorem}
This yields immediately the bound
\begin{equation}
\label{eqn:complex}
N_{\C}(s,d,k) \leq \sum_{0 \leq \ell \leq k} {s \choose k-\ell} d^{k - \ell}.
\end{equation}

Notice that
\[
\sum_{0 \leq \ell \leq k} {s \choose k-\ell} d^{k - \ell}
\leq {sd + k \choose k}
\]
for all values of $s,d,k \geq 0$, and thus
the bound in (\ref{eqn:complex}) is slightly better than that in 
(\ref{eqn:rbg}).

\begin{remark}
We would like to point out that 
it is tempting to try to prove Theorem \ref{the:main1} 
by first proving the bounds for sufficiently generic families of 
polynomials ${\mathcal P} \subset \R[X_1,\ldots,X_k]$,
and then proving that generic families of polynomials actually 
represents the worst case. While it is indeed easy to prove the bound
in  Theorem \ref{the:main1} for generic families,
it is not at all clear how to prove the second statement.

In some special cases, for instance for families of linear polynomials,
one can prove that infinitesimal perturbations of the input polynomials
can only increase the number of cells in the arrangement
(see for instance \cite{Matousek}). Hence, for the problem of bounding the
number of cells in an arrangement of hyperplanes it
suffices to consider hyperplanes in general position. 
But for polynomials of degree greater than one, 
there is no guarantee that 
an arbitrary infinitesimal perturbation of any given family of 
polynomials will not 
decrease
the Betti numbers of the realizations
of the realizable sign conditions.
\end{remark}

The rest of the paper is organized as follows. In Section \ref{sec:prelim},
we state some results we will need for the proof of the main theorems with
appropriate references. In Section \ref{sec:main} we prove the main results
of the paper.

\section{Mathematical Background}
\label{sec:prelim}
In this section we state without proof some well known results regarding the
topology of complex affine varieties which we will need in the 
proof of the main theorem. We also point the reader to appropriate
references for the proofs of these results.

\subsection{Homotopy type of complex affine varieties}
Let $V \subset \C^k$ be an affine variety of complex dimension $\ell$. 
If we consider $V$ as a real semi-algebraic set in $\R^{2k}$, then its
real dimension is $2\ell$. 
It was proved by Karchyauskas \cite{Kr} that in the case $\C = \complex$,
$V$ has the homotopy type of a CW-complex of (real) dimension $\ell$. 
A proof which uses stratified Morse
theory can be found in  \cite[Section 5.1, page 198]{GM}. 

As a consequence of above results we have that the 
homology groups of $V$, $\HH_i(V)$, vanish for all $i > \ell$. 
By a standard
transfer argument using the Tarski-Seidenberg transfer principle
this latter statement
can be extended easily to all algebraically closed fields of char $0$ 
(see for instance \cite[pp. 261]{BPRbook2} for an argument of this type).
As a result we obtain the following theorem.

\begin{theorem}
\label{the:CW}
Let $V \subset \C^k$ be an affine variety of (complex) dimension $\ell$.
Then $\HH_i(V) = 0$ for all $i > \ell$.
\end{theorem}

\begin{remark}
Note that in the case  $\C = \complex$, 
Theorem \ref{the:CW} holds much more generally for 
{\em Stein spaces}  (see \cite{Narasimhan67}). However, we do not need to use
this fact.
\end{remark}

\subsection{Generalized Bezout Inequality}
In our proof we will need to use the following 
generalized form of Bezout's theorem
(see for instance \cite[Example 8.4.6, page 148]{Fulton}).

Recall that the degree of an irreducible complex projective variety 
$V \subset \PP_{\C}^k$ of complex dimension $\ell$, denoted
$\deg(V)$,  is the number
of points in the intersection $V \cap L$, where $L$ is a generic
linear subspace $\PP_{\C}^k$ of dimension $k-\ell$.

If $V = \cup_{1 \leq i \leq n} V_i$ is a union of irreducible varieties
each of dimension $\ell$, 
then we define $\deg(V) = \sum_{1 \leq i \leq n}
\deg(V_i)$.

\begin{theorem}[Generalized Bezout]
\label{the:bezout}
Let $V_1,\ldots,V_s$ be complex projective varieties,
such that each $V_i$ is a union of irreducible varieties of the same dimension,
and 
$Z_1,\ldots,Z_t$ be the irreducible components of the intersection
$V_1\cap \cdots \cap V_s$. Then
$$
\displaylines{
\sum_{1 \leq i \leq t} \deg(Z_i) \leq \prod_{1 \leq i \leq s} \deg(V_i).
}
$$
\end{theorem}

\subsection{Exact Sequence for Borel-Moore Homology}
We will assume that the reader is familiar with the definition of 
Borel-Moore homology groups \cite{BM} for locally closed semi-algebraic sets
(see for instance \cite[Definition 11.7.13]{BCR}).
Let $W \subset \C^k$ be a closed set and $A \subset W$ a 
closed subset of $W$. Then the following exact sequence
relates  the Borel-Moore homology groups of $W$, $A$ and $W\setminus A$
(see \cite[Proposition 11.7.15]{BCR}).
\[
\cdots \rightarrow \HH_p^{BM}(A) \rightarrow \HH_p^{BM}(W) \rightarrow 
\HH_p^{BM}(W\setminus A)
\rightarrow \HH_{p-1}^{BM}(A) \rightarrow \cdots
\]
It follows that
\begin{proposition}
\label{prop:BM}
For each $p > 0$
\[
b_p^{BM}(W \setminus A) \leq b_{p-1}^{BM}(A) + b_p^{BM}(W).
\]
In particular,
\[
\sum_{p \geq 0} b_p^{BM}(W \setminus A) \leq \sum_{p \geq 0} b_{p}^{BM}(A) 
+ \sum_{p \geq 0} b_p^{BM}(W).
\]
\end{proposition}

One consequence of Proposition \ref{prop:BM} (and Theorem \ref{the:CW}) is
the following important fact.

\begin{proposition}
\label{prop:affine}
Let $\rho \in \mbox{Zero-pattern}({\mathcal P})$ and let
$d_\rho$ be the complex dimension of $\RR(\rho,\C^k)$. Then
$\displaystyle{
b_i^{BM}(\rho,\C^k)  = 0 \mbox{ for all } i > d_\rho.
}
$
\end{proposition}

\begin{proof}
Let
$$
\displaylines{
W' = \bigcap_{P\in{\mathcal P},\rho(P) =0} \ZZ(P,\C^k),
}
$$
and 

$$
\displaylines{
A'  = W' \cap 
\left(\bigcup_{P\in{\mathcal P},\rho(P) =1}\ZZ(P,\C^k)\right).
}
$$
It is clear that $W'$ (resp. $A'$)
is an  affine sub-variety  of $\C^k$ of complex dimension
$d_\rho$ (resp. $ < d_\rho$).

Now by definition (cf. Eqn. \eqref{eqn:C(Z)})
\begin{eqnarray*}
\RR(\rho,\C^k)   &=&      \{x\in \C^k\;\mid\; 
\bigwedge_{P\in{\mathcal P}} \sz({P}(x))=\rho(P) \}
\cap \overline{B_{2k}(0,\Omega)} \\
&=& 
\left(W' \cap \overline{B_{2k}(0,\Omega)}
\right) \setminus 
\left(A'\cap \overline{B_{2k}(0,\Omega)} \right)
.
\end{eqnarray*}

Now applying Proposition \ref{prop:BM} with
\[
W = W' \cap \overline{B_{2k}(0,\Omega)},
\]
and 
\[
A = A' \cap \overline{B_{2k}(0,\Omega)}
\]
and observing that by Theorem \ref{the:CW}
\[
b_i^{BM}(A) = b_i(A) = 0 \mbox{ for all }i > d_\rho -1,
\] 
and
\[
b_i^{BM}(W) = b_i(W) = 0 \mbox{ for all } i > d_\rho,
\] 
we finally obtain that 
\[
b_i^{BM}(\rho,\C^k)  = b_i^{BM}(W\setminus A) = 0 \mbox{ for all } i > d_\rho.
\] 
\end{proof}

Now suppose that $\displaystyle{A = \bigcup_{1 \leq i \leq n} A_i}$, 
where each $A_i$ is a closed and bounded 
semi-algebraic subset of $\C^k = \R^{2k}$. 
The following inequality is
a consequence of the Mayer-Vietoris exact sequence for simplicial
homology groups (noting that the Borel-Moore homology are isomorphic
to the simplicial homology groups for closed and bounded 
semi-algebraic
sets).

\begin{proposition}
\label{prop:MV}
For each $p \geq 0$
\[
b_p^{BM}(A) \leq \sum_{0 \leq i \leq p} \sum_{1 \leq \alpha_0 < 
\cdots <\alpha_i \leq n}b_{p-i}^{BM}(A_{\alpha_0}\cap\cdots\cap A_{\alpha_i}).
\]
\end{proposition}

\subsection{Smith Inequalities}
Let $X$ be any compact topological space and $c: X \rightarrow X$ an
involution. Let $F \subset X$ denote the set of fixed points of $c$,
and let $X'$ denote the quotient space 
$X/(x \sim  cx)$. 
The projection map
$X \rightarrow X'$ maps $F$ homeomorphically onto a subset of $X'$ 
(which we also denote by $F$). The following exact sequence called the Smith
exact sequence is well known (see for instance \cite[page 131]{Viro}).
\[
\cdots \rightarrow \HH_p(X',F)\oplus \HH_p(F) \rightarrow \HH_p(X) \rightarrow \HH_p(X',F)\rightarrow \HH_{p-1}(X',F)\oplus \HH_{p-1}(F) \rightarrow \cdots
\]
(here as elsewhere in this paper all homology groups are taken with 
coefficients in ${\mathbb Z}_2$).
As an immediate  consequence we have that for any $q \geq 0$
\begin{equation}
\label{eqn:smith}
\sum_{i \geq q} b_i(F) \leq \sum_{i \geq q} b_i(X).
\end{equation}

We are going to apply 
Inequality~\eqref{eqn:smith}
in the case where 
$X \subset \C^k$ 
is the intersection with $\overline{B_{2k}(0,\Omega)}$ of the 
basic constructible subset defined by a formula
\begin{equation*}
\label{eqn:formula}
P_1 = \cdots P_\ell = 0, P_{\ell+1}\neq 0,\ldots, P_{s} \neq 0,
\end{equation*}
where $P_i \in \R[X_1,\ldots,X_k]$.
Let $c$ denote the complex conjugation and $F \subset X$ its set of
fixed points. 
Clearly, 
$F \subset \overline{B_k(0,\Omega)}$  is defined
by the same formula as in (\ref{eqn:formula}).

The Borel-Moore Betti numbers of $X$ (as well as $F$) are by definition
the ordinary simplicial homology groups of the compact pairs
$(\overline{X},\overline{X}\setminus X)$ (resp.
$(\overline{F},\overline{F}\setminus F)$) where
$\overline{X}$ (resp. $\overline{F}$) denotes the closure 
(in the Euclidean topology) of
$X$ (respectively $F$).
It is easy to see that the action of  conjugation $c$ 
extends to $\overline{X}$ and  the Smith exact sequence yields the 
following inequality.

\begin{proposition}[Smith Inequality]
\label{prop:smith}
For any $q \geq 0$
\[ 
\sum_{i \geq q} b_i^{BM}(F) \leq \sum_{i \geq q} b_i^{BM}(X).
\]
\end{proposition}
   
\subsection{Oleinik-Petrovsky-Thom-Milnor Inequalities}
We will also need the following inequality proved separately by 
Oleinik and Petrovsky \cite{OP}, Thom \cite{T} and Milnor \cite{Milnor2},
bounding the sum of the Betti numbers of a real algebraic set. In particular,
it also gives a bound on the sum of the Betti numbers of a complex algebraic
set considered as a real algebraic set in an affine space twice the dimension
of the complex one.
Moreover, since Borel-Moore homology agrees with singular homology for
closed and bounded
sets, we have 
\begin{proposition}
\label{prop:OP1}
Let ${\mathcal Q} \subset
 \R[X_1,\ldots,X_k],$ $\deg(Q) \leq  d, Q \in {\mathcal Q}$,
$V_{\R} = \ZZ({\mathcal Q},\R^k) \cap \overline{B_k(0,\Omega)}$,
and $V_{\C} = \ZZ({\mathcal Q},\C^k) \cap \overline{B_{2k}(0,\Omega)}$. Then
\begin{eqnarray}
\sum_{0 \leq i \leq k} b_i(V_{\R}) &=& \sum_{0 \leq i \leq k} b_i^{BM}(V_{\R}) 
\leq d(2d-1)^{k-1} \\
\sum_{0 \leq i \leq k} b_i(V_{\C}) &=&
\sum_{0 \leq i \leq k} b_i^{BM}(V_{\C}) 
\leq d(2d-1)^{2k-1}.
\end{eqnarray}
\end{proposition}

\section{Proofs of the Main Results}
\label{sec:main}
As mentioned before, we first prove a bound on the sum of Borel-Moore
Betti numbers of the realizations of all realizable zero patterns of
a given family of polynomials over $\C^k$. We will then apply this
result in the real case via the Smith inequality (Proposition \ref{prop:smith}). 
The first part of 
this section is devoted to the proof of the following proposition.

\begin{proposition}
\label{prop:complexmain}
Let ${\mathcal P} \subset \C[X_1,\ldots,X_k]$ be a family of polynomials
with $\#{\mathcal P} = s$ and $\deg(P) \leq d$ for all $P \in {\mathcal P}$.
Then

\begin{eqnarray*}
\sum_{\rho}
b^{BM}(\rho,\C^k)  & \leq &  
\sum_{0 \leq \ell \leq k} \left(
{s \choose k - \ell}{s  \choose \ell}d^k 
 + \sum_{1 \leq i \leq \ell} {s \choose k - \ell}{s \choose \ell - i}d^{O(k^2)}
\right)
\\
&=& 
\sum_{0 \leq \ell \leq k} {s \choose k - \ell}{s  \choose \ell}d^k + 
O(s^{k-1})
\end{eqnarray*}
where the sum on the left hand side is taken over all $\rho \in  
\mbox{Zero-pattern}({\mathcal P})$.
\end{proposition}

For $\rho \in \mbox{Zero-pattern}({\mathcal P})$
let ${\mathcal P}_\rho = \{P \in {\mathcal P} \mid \rho(P) = 0\}$ and
let $V_\rho = \ZZ({\mathcal P}_\rho, \C^k).$ Let 
$d_\rho = \dim_{\C} \RR(\rho,\C^k)$. 
Let $V_\rho = V_{\rho,1} \cup \cdots \cup V_{\rho,n(\rho)}$ denote the
decomposition of $V_\rho$ into irreducible components.

We first prove a bound on the number of isolated points occurring 
as realizations of zero patterns.
\begin{lemma}
\label{lem:isolatedpoints}
\[
\# \{(\rho,i) \mid \dim_{\C} V_{\rho,i} = 0\} \leq {s \choose k} d^k.
\]
\end{lemma}

\begin{proof}
Let ${\mathcal P} = \{P_1,\ldots,P_s\}$ (the implicit ordering of
the polynomials in ${\mathcal P}$ will play a role in the proof).

First consider a sequence of $k$ polynomials of ${\mathcal P}$, 
$P_{i_1},\ldots,P_{i_k},$ with
$1 \leq i_1 \leq \cdots \leq i_k \leq s$.
Let $\{W_{\alpha_1}\}_{\alpha_1 = 1,2..}$ denote the irreducible
components of $\ZZ(P_{i_1},\C^k)$. Similarly, let
$\{W_{\alpha_1,\alpha_2}\}_{\alpha_2 = 1,2...}$ denote the 
irreducible components  of $W_{\alpha_1} \cap \ZZ(P_{i_2},\C^k)$ and so on.
We prove by induction that for $1 \leq \ell \leq k$
\begin{equation}
\label{eqn:bound}
\sum_{\alpha_1,\alpha_2,\ldots,\alpha_\ell} \deg(W_{\alpha_1,\ldots,\alpha_\ell}) \leq d^\ell.
\end{equation}

The claim is trivially true for $\ell=1$. Suppose it is true up to $\ell -1 $.
Then using Theorem \ref{the:bezout} we get that
\[
\sum_{\alpha_1,\ldots,\alpha_{\ell-1}}
 {\deg(W_{\alpha_1,\ldots,\alpha_{\ell-1}} \cap \ZZ(P_\ell,\C^k))}  
 \leq d^{\ell-1}\cdot d = d^\ell
\]

Now let $\dim_{\C}V_{\rho,i} = 0$ and we denote the corresponding point in
$\C^k$ by $p$.

We first claim that
there must exist $P_{i_1},\ldots,P_{i_k} \in {\mathcal P}_{\rho}$, with
$1 \leq i_1 < i_2 < \cdots < i_k \leq s$,  and
irreducible affine algebraic varieties 
$V_0, V_1,\ldots,V_k$ satisfying:
\begin{enumerate}
\item
$\C^k= V_0 \supset V_1 \supset V_2 \supset \cdots \supset V_k = \{p\}$, and
\item
for each $j, \;1 \leq j \leq k$,
$V_{j}$ is an irreducible component of $V_{j-1} \cap \ZZ(P_{i_j},\C^k)$
which contains $p$.
\end{enumerate} 
Note that in the sequence, $V_0,V_1,\ldots,V_k$,
$\dim_{\C}(V_i)$ is necessarily  equal to $k-i$.
To prove the claim, let 
\[
{\mathcal P}_\rho = \{P_{j_1},\ldots,P_{j_m}\},
\]
with $1 \leq j_1 < j_2 < \cdots < j_m \leq s$.
Let $W_0 = \C^k$ and
for each $h, 1 \leq  h \leq m$ define inductively
$W_{h}$ to be an irreducible component of $W_{h-1} \cap \ZZ(P_{j_h},\C^k)$
containing the point $p$. By definition,
$\C^k = W_0 \supset W_1 \supset \cdots \supset W_m = \{p\}$ and each
$W_h$ is irreducible. Hence, there must exist precisely
$k$ indices,
$1 \leq \alpha_1 < \cdots < \alpha_k \leq m$ such that
$\dim_{\C}W_{\alpha_{i}} = \dim_{\C} W_{\alpha_{i-1}} - 1$.
If $\beta \not\in \{\alpha_1,\ldots,\alpha_k\}$, then
$W_{\beta} = W_{\beta - 1}$.
Now let $i_1 = j_{\alpha_1},\ldots, i_k = j_{\alpha_k}$,
and $V_j = W_{i_j}$ for $0 \leq j \leq k$. It is clear from construction
that the sequence $V_0,\ldots,V_k$ satisfies the properties stated above,
thus proving the claim.
 
The number of sequences of
polynomials of length $k$,
$P_{i_1},\ldots,P_{i_k},$ with
$1 \leq i_1 < \cdots < i_k \leq s$,  is ${s \choose k}$. The lemma
now follows from the bound in (\ref{eqn:bound}).
\end{proof} 

The following lemma 
%%sb adds
(which also appears in \cite{JS00})
is a generalization of the above lemma to positive
dimensional realizations of zero patterns.

\begin{lemma}
\label{lem:positive}
For $0 \leq \ell \leq k$
$$
\displaylines{
\sum_{\rho \in \mbox{\rm Zero-pattern}({\mathcal P})} 
\sum_{1 \leq i \leq n(\rho), \dim_{\C}V_{\rho,i} =\ell} 
\deg V_{\rho,i} \leq {s \choose k-\ell} d^{k - \ell}.
}
$$
\end{lemma}

\begin{proof}
In order to bound the degree of $V_{\rho,i}$ with $\dim_{\C}V_{\rho,i} =\ell$
it suffices to count the isolated zeros of the
intersection of $V_{\rho,i}$ with a generic affine subspace $L$ of dimension
$k-\ell$. Since there are only a finite number of $V_{\rho,i}$'s to consider
we can assume that
\begin{enumerate}
\item  $L$ is generic for all of them,
\item
the sets $V_{\rho,i} \cap L$ are pairwise disjoint,
\item
for each $V_{\rho,i}$ with $\dim_{\C}V_{\rho,i} =\ell$,
$L \cap V_{\rho,i} \cap \ZZ(P,C^k) = \emptyset $ for each 
$P \in {\mathcal P} \setminus {\mathcal P}_\rho.$
\end{enumerate}
Now restrict to the subspace $L$ and apply Lemma 
\ref{lem:isolatedpoints},
noting that an isolated point of
$V_{\rho,i} \cap L$ is an isolated point of the family
${\mathcal P}$ restricted to $L$.
\end{proof}

For each $\rho \in \mbox{Zero-pattern}({\mathcal P})$, let
$W_{\rho,1},\ldots,W_{\rho,m(\rho)}$ be the irreducible components
of $V_\rho$ of dimension $d_\rho$ having the property that no polynomial in
${\mathcal P} \setminus {\mathcal P}_\rho$ vanishes identically 
on any of them. This implies that each $W_{\rho,j} \cap \RR(\rho,\C^k)$ is
non-empty and of complex dimension $d_\rho$. 
(Thus, the union of the sets $W_{\rho,j} \cap \RR(\rho,\C^k)$ is the full
dimensional part of the constructible set $\RR(\rho,\C^k)$.)
Let 
$$
\displaylines{
W_\rho = \bigcup_{1 \leq j \leq m(\rho)} W_{\rho,j}.
}
$$

We first note that 
$b_\ell^{BM}(W_\rho)$ 
can be bounded in terms of 
$d$ and $k$ independent of $s$ as follows.
The exact nature of the dependence of this bound on
$d$ and $k$ is not important for the asymptotic result that we prove 
in this paper. However, we are able to show the following.

\begin{lemma}
\label{lem:OP1}
$$
\displaylines{
\sum_{i \geq 0} 
b_i^{BM}(W_\rho \cap \overline{B_{2k}(0,\Omega)}) \leq d^{O(k^2)}.
}
$$
\end{lemma}

\begin{proof}
We first claim that $\deg(W_\rho) \leq d^{O(k)}.$ To see this first 
observe that $\deg(W_\rho) \leq \deg(V_\rho')$ where
$
\displaystyle{
V_\rho' = \bigcup_{1 \leq i \leq n(\rho), \dim_{\C}V_{\rho,i} = d_\rho} 
V_{\rho,i}.
}$
Now $\deg(V_{\rho}')$ is the number of isolated points in the intersection
of $V_{\rho}$ with a generic affine subspace, $L \subset \C^k$, 
of dimension $k-d_\rho.$
Now $V_\rho$ is an affine algebraic set in $\C^k$ defined by polynomials
of degree at most $d$. Identifying $\C^k$ with $\R^{2k}$, and
separating real and imaginary parts, we have that the intersection 
$L \cap V_{\rho} \subset \R^{2(k - d_\rho)}$, is a real algebraic
set defined by polynomials of degree at most $d$. Applying 
Proposition \ref{prop:OP1} we have that the number of isolated points of
 $L \cap V_{\rho}$ is at most $b_0(L \cap V_{\rho}) \leq d(2d - 1)^{2k-1}$.
Hence, $\deg(W_\rho) \leq d^{O(k)}.$

Since $W_\rho$ is an affine variety of
degree at most $d^{O(k)}$, it is defined by at most $k+1$ polynomials each of
degree at most $d^{O(k)}$. 

Now using Proposition \ref{prop:OP1} we obtain that
\[
\sum_{i} b_i^{BM}(W_{\rho} \cap \overline{B_{2k}(0,\Omega)}) \leq (d^k)^{O(k)} = d^{O(k^2)}.
\] 
\end{proof}

If ${\mathcal Q}$ is any subset of ${\mathcal P}$ an easy extension of the
above argument also yields
\begin{lemma}
\label{lem:OP2}
$$
\displaylines{
\sum_{i \geq 0}
b_i^{BM}(W_\rho \cap \ZZ({\mathcal Q},\C^k) \cap \overline{B_{2k}(0,\Omega)})
\leq d^{O(k^2)}.
}
$$
\end{lemma}

We remark here that the bounds in \ref{lem:OP1} and \ref{lem:OP2} are unlikely
to be tight and improving them might lead to bounds which are 
optimal not just asymptotically, but for all values of $s,d$ and $k$. 

The following lemma is key to the proof Theorem \ref{the:main1}.

\begin{lemma}
\label{lem:onecondition}
Let $\ell >0$ and $\rho \in \mbox{Zero-pattern}({\mathcal P})$ with 
$d_\rho = \ell$. Then
\begin{eqnarray*}
b^{BM}_\ell(\rho,\C^k) &\leq & \sum_{1 \leq i \leq m(\rho)}
\deg(W_{\rho,i}){s \choose \ell} d^\ell + 
\sum_{1 \leq i \leq \ell}
{s \choose \ell-i } d^{O(k^2)} \\
&=& \sum_{1 \leq i \leq m(\rho)}
\deg(W_{\rho,i}){s \choose \ell} d^\ell + 
O(s^{\ell - 1}).
\end{eqnarray*}

For $j < \ell$
\begin{eqnarray*}
b^{BM}_j(\rho,\C^k) &\leq & 
\sum_{1 \leq i \leq \ell}
{s \choose \ell-i } d^{O(k^2)} \\
&=& O(s^{\ell-1}).
\end{eqnarray*}
\end{lemma}

\begin{proof}
Let
$\rho \in \mbox{Zero-pattern}({\mathcal P})$ with $ d_\rho = \ell,$
and let $A_\rho = (W_{\rho} \cap \overline{B_{2k}(0,\Omega)}) \setminus \RR(\rho,\C^k).$ 
Clearly
\[
A_\rho = (W_\rho \cap \overline{B_{2k}(0,\Omega)}) \cap (\cup_{P \in {\mathcal P}\setminus {\mathcal P}_\rho}
\ZZ(P,\C^k)). 
\]

Now
\[
b_\ell^{BM}(\rho,\C^k) = b_\ell^{BM}(W_\rho \cap \RR(\rho,\C^k))
\]
by Proposition \ref{prop:affine}
since $W_\rho \cap \RR(\rho,\C^k)$
is the top dimensional part of $\RR(\rho,\C^k)$.

Also notice that
\[
W_\rho \cap \RR(\rho,\C^k)= (W_\rho \cap \overline{B_{2k}(0,\Omega)}) \setminus A_\rho.
\]

Moreover, from Proposition \ref{prop:BM} it follows that
\[
b_\ell^{BM}((W_\rho \cap \overline{B_{2k}(0,\Omega)}) \setminus A_\rho) \leq 
b_\ell^{BM}(W_\rho \cap \overline{B_{2k}(0,\Omega)}) + 
b_{\ell - 1}^{BM}(A_\rho).
\]

Now observe that by Lemma \ref{lem:OP1}
\[
b_\ell^{BM}(W_\rho \cap \overline{B_{2k}(0,\Omega)}) \leq d^{O(k^2)}.
\]

We now   bound the second term
$ b^{BM}_{\ell - 1}(A_\rho)=b_{\ell - 1}(A_\rho)$ 
(since $A_\rho$ is closed and bounded) in terms of $s$ as follows.

Expressing $A_\rho$  as the union
$
\displaystyle{
\bigcup_{P \in {\mathcal P}\setminus {\mathcal P}_\rho}
\ZZ(P,\C^k) \cap W_\rho \cap \overline{B_{2k}(0,\Omega)}
}$,
we obtain using Mayer-Vietoris inequalities 
(Proposition \ref{prop:MV})
$$
\displaylines{
b_{\ell -1}(A_\rho) \leq \sum_{ {\mathcal Q}
\subset {\mathcal P}\setminus {\mathcal P}_\rho,\#{\mathcal Q} = \ell}
b_0(W_\rho \cap \ZZ({\mathcal Q},\C^k) \cap \overline{B_{2k}(0,\Omega)})\cr  
+ \sum_{1 \leq i \leq \ell}
\sum_{ {\mathcal Q}
\subset {\mathcal P}\setminus {\mathcal P}_\rho,\#{\mathcal Q} = \ell-i}
b_i(W_\rho \cap \ZZ({\mathcal Q},\C^k)\cap \overline{B_{2k}(0,\Omega)}). 
}
$$

Now using the generalized Bezout inequality (Theorem \ref{the:bezout}) and
noting that for any complex affine algebraic set  $X \subset \C^k$,
$b_0(X) \leq \deg(X)$, we have that
\[
b_0(W_\rho \cap \ZZ({\mathcal Q},\C^k)) 
= b_0(W_\rho \cap \ZZ({\mathcal Q},\C^k) \cap \overline{B_{2k}(0,\Omega)}) 
\leq 
\deg(W_\rho) d^\ell,
\]
since each polynomial in ${\mathcal Q}$ has  degree at most $d$.
Thus, 
\[
\sum_{ {\mathcal Q}
\subset {\mathcal P}\setminus {\mathcal P}_\rho,\#{\mathcal Q} = \ell}
b_0(W_\rho \cap \ZZ({\mathcal Q},\C^k) \cap \overline{B_{2k}(0,\Omega)})
\leq 
\deg(W_\rho) {s  \choose \ell} d^\ell.
\]
Also, using Lemma \ref{lem:OP2} we have that the second part of the sum
\[
\sum_{1 \leq i \leq \ell}
\sum_{ {\mathcal Q}
\subset {\mathcal P}\setminus {\mathcal P}_\rho,\#{\mathcal Q} = \ell-i}
b_i(W_\rho \cap \ZZ({\mathcal Q},\C^k) \cap \overline{B_{2k}(0,\Omega)}) 
\leq \sum_{1 \leq i \leq \ell} {s \choose \ell - i} d^{O(k^2)}.
\]
\end{proof}

\begin{lemma}
\label{lem:complex}
Let $0 < \ell \leq k$. 
Then
\begin{eqnarray*}
\sum_{\rho}
b^{BM}(\rho,\C^k)  & \leq & 
{s \choose k - \ell}{s  \choose \ell}d^k + 
\sum_{1 \leq i \leq \ell} {s \choose k - \ell}{s \choose \ell - i}d^{O(k^2)}
\\
&=& 
{s \choose k - \ell}{s  \choose \ell}d^k + 
O(s^{k-1})
\end{eqnarray*}
where the sum on the left hand side is taken over all 
$\rho  \in \mbox{Zero-pattern}({\mathcal P})$ with $d_\rho = \ell $.
\end{lemma}

\begin{proof}
Follows from Lemmas \ref{lem:onecondition} and \ref{lem:positive}
and the fact that
$b_i^{BM}(\rho,\C^k) = 0$ for $i > \ell$ for all 
$\rho  \in \mbox{Zero-pattern}({\mathcal P})$ with $d_\rho = \ell $
(Proposition \ref{prop:affine}).
\end{proof}

\begin{proof}[Proof of Proposition \ref{prop:complexmain}]
Follows directly from 
 from Lemma \ref{lem:complex},
since 
$0 \leq d_\rho \leq k$, for each $\rho \in \mbox{Zero-pattern}({\mathcal P})$.
\end{proof}

We now consider the real case directly. We first prove a few preliminary 
lemmas.

\begin{lemma}
\label{lem:1}
Let $C \subset \R^k$ be a bounded semi-algebraically connected component having real dimension
$\ell$ of the basic  semi-algebraic set defined by
$$
\displaylines{
P_1 = \ldots = P_m= 0, \cr
P_{m+1} > 0,\ldots, P_n > 0, \cr
P_{n+1} < 0,\ldots,P_s < 0.
}
$$
Let $h: \Delta \rightarrow C$ be a semi-algebraic triangulation of $C$ and
$\sigma$ an $(\ell-1)$-simplex of the triangulation. Then the number of 
$\ell$-simplices $\eta$ in $\Delta$ such that $\sigma$ is a face of
$\eta$ is even.
\end{lemma}   
\begin{proof}
Consider the real algebraic set $V \subset \R^{k+s-m}$ defined by
$$
\displaylines{
P_1 = \cdots = P_m = 0,\cr
T_{m+1}^2P_{m+1} - 1= \cdots = T_{n}^2P_n - 1 = 0, \cr
T_{n+1}^2P_n + 1 = \cdots = T_s^2 P_s +1 = 0,
}
$$
where $T_{m+1},\ldots,T_s$ are new variables. Let
$\pi: \R^{k+s-m} \rightarrow \R^k$ be the projection map which forgets
the new co-ordinates.
Since $C$ is bounded, $\pi^{-1}(C) \cap V$ is the disjoint union of 
$2^{s-m}$  semi-algebraically connected components, $C_1,\ldots,C_{2^{s-m}}$ of $V$, 
each homeomorphic to $C$. Moreover, any 
triangulation of $C$ can be lifted to a triangulation of 
$\cup_{1 \leq i \leq 2^{s-m}} C_i$.

Now, $C_1$ is a semi-algebraically connected component of a real algebraic set of dimension 
$\ell$. Thus, it suffices to prove the lemma in case $C$ is a semi-algebraically connected
component of a real algebraic set. We now refer the reader to the proof
of Theorem 11.1.1 in \cite{BCR}, where the same claim is proved in case 
$C$ is a real algebraic set of dimension $\ell$. However, the proof which 
is of a local character
also applies in case $C$ is only a semi-algebraically connected component of a real algebraic set.
\end{proof}

\begin{lemma}
\label{lem:2}
Let $C \subset \R^k$ be a bounded semi-algebraically connected component of a basic
semi-algebraic set, $S \subset \R^k$,  defined by
$P_1 = \cdots = P_m = 0, P_{m+1} > 0,\ldots, P_{n} > 0, P_{n+1} < 0,
\ldots, P_{s} < 0$, such that the real dimension of $C$ is  $\ell$. 
Then $b_\ell^{BM}(C) \geq 1$. 
\end{lemma}

\begin{proof}
Let $S \subset \R^k$ be the set defined by 
$P_1 = \cdots = P_m = 0, P_{m+1} > 0,\ldots, P_{n} > 0, P_{n+1} < 0,
\ldots, P_{s} < 0$ and let $T = S \cap \overline{B_k(0,\Omega)}$. 
For all sufficiently large $\Omega > 0$
each bounded semi-algebraically connected component of $S$,
is also a semi-algebraically connected component of $T$.
Since $T$ is bounded, from the definition of Borel-Moore homology groups
we have that
\[
\HH_i^{BM}(T) \cong \HH_i(\overline{T}, \overline{T} \setminus T).
\]
Let $h: \Delta \rightarrow \overline{T}$ be a semi-algebraic triangulation
of the pair $(\overline{T}, \overline{T} \setminus T)$.
 
Now by Lemma \ref{lem:1} 
each $(\ell -1)$-simplex $\sigma$ in $\Delta$ whose image is 
contained in $C$ is a face of an 
even number $\ell$-simplices 
whose images are contained in
in $C$. Moreover, 
if $\sigma$ is a face of an $\ell$-simplex $\eta$ whose image is in 
$C$ and $h(\sigma) \not\subset C$ then 
$h(\sigma) \subset \overline{T} \setminus T$.
It now follows that there exists a non-empty set of $\ell$-simplices
whose images are 
contained in $C$, having the property that 
each $(\ell-1)$-face of a simplex in the family is either contained
in an even number of simplices of the family or has its image contained in 
$\overline{T} \setminus T$. Let $Z \in C_\ell(\Delta)$ 
denote the sum of these simplices. Clearly, $Z$ represents a 
non-zero $\Z_2$-homology class in $\HH_\ell^{BM}(C)$.
Thus, 
\[
b_\ell^{BM}(C) \geq 1.
\]
\end{proof}

Let ${\mathcal P} \subset \R[X_1,\ldots,X_k]$. For 
${\mathcal Q} \subset {\mathcal P}$ let 
\[
 \Sigma_{\mathcal Q} = \{
\sigma \in \{0,1,-1\}^{\mathcal P} \;\mid\; 
\sigma(P) = 0 \;\mbox{ for } P \in {\mathcal Q}, \mbox{ and }\;\sigma(P) \neq 0 \; \mbox{ for }  P \in
{\mathcal P} \setminus {\mathcal Q}\},
\]
and let $\rho_{\mathcal Q} \in \{0,1\}^{\mathcal P}$
be defined by
\begin{eqnarray*} 
\rho_{\mathcal Q}(P) &=&  0,\; \mbox{ for }  P \in {\mathcal Q},  \\
&=& 1, \; \mbox{ for }  P \in {\mathcal P} \setminus {\mathcal Q}.
\end{eqnarray*}

Also, let 
$$
\displaylines{
W_{{\mathcal Q},\R} =
\bigcup_{\sigma \in  \Sigma_{\mathcal Q}} \RR(\sigma,\R^k).
}
$$

\begin{lemma}
\label{lem:3}
If $C$ is a non-empty semi-algebraically connected component of
$\RR(\sigma,\R^k)$
for some $\sigma \in \Sigma_{\mathcal Q},$
and 
$i :C \hookrightarrow W_{{\mathcal Q},\R}$ the inclusion map,
then
$i_*: \HH_i^{BM}(C) \rightarrow \HH_i^{BM}(W_{{\mathcal Q},\R})$
is an injection if $i > 0$ and
in case $i=0$, $i_*$ is injective if $C$ is closed.
In particular, $i_*: \HH_0^{BM}(C) \rightarrow \HH_0^{BM}(W_{{\mathcal Q},\R})$ 
is injective if $C$ is zero dimensional.
\end{lemma}

\begin{proof}
Follows directly from the definition of Borel-Moore homology groups.
\end{proof}

\begin{proof}[Proof of Theorem \ref{the:main1}]
We follow the notation introduced above.
Let ${\mathcal Q} \subset {\mathcal P} \subset \R[X_1,\ldots,X_k]$. 

Notice that $W_{{\mathcal Q},\R}$ is the fixed point of set of
the action of complex conjugation on $\RR(\rho_{\mathcal Q},\C^k)$.
It follows from Proposition \ref{prop:smith} that

\begin{equation}
\label{eqn:proofmain3}
b^{BM}(W_{{\mathcal Q},\R}) 
\leq
b^{BM}(\rho_{\mathcal Q},\C^k).
\end{equation}

Applying Proposition \ref{prop:complexmain} we obtain that
\begin{eqnarray*}
\label{eqn:proofmain4}
\sum_{\rho}
b^{BM}(\rho,\C^k)  & \leq & 
\sum_{0 \leq \ell \leq k}
\left({s \choose k - \ell}{s  \choose \ell}d^k + 
\sum_{1 \leq i \leq \ell} {s \choose k - \ell}{s \choose \ell - i}d^{O(k^2)}
\right)\\
&=& 
\sum_{0 \leq \ell \leq k} {s \choose k - \ell}{s  \choose \ell}d^k + 
O(s^{k-1})
\end{eqnarray*}
where the sum on the left hand side is taken over all 
$\rho  \in \mbox{Zero-pattern}({\mathcal P})$.

Now applying Inequality \eqref{eqn:proofmain3} we get
\begin{eqnarray}
\label{eqn:proofmain5}
\sum_{{\mathcal Q} \subset {\mathcal P}} b^{BM}(W_{{\mathcal Q},\R})
& \leq &  
\sum_{{\mathcal Q} \subset {\mathcal P}} b^{BM}(\rho_{\mathcal Q},\C^k)
\nonumber\\
& = &
\sum_{\rho \in \mbox{Zero-pattern}({\mathcal P})} b^{BM}(\rho,\C^k) \nonumber\\
& \leq & 
\sum_{0 \leq \ell \leq k}\left({s \choose k - \ell}{s  \choose \ell}d^k + 
\sum_{1 \leq i \leq \ell} {s \choose k - \ell}{s \choose \ell - i}d^{O(k^2)}
\right)\nonumber\\
& \leq &
\sum_{0 \leq \ell \leq k} {s \choose k - \ell}{s  \choose \ell}d^k + 
O(s^{k-1}).
\end{eqnarray}

We now use Inequality \eqref{eqn:proofmain5} to 
bound the number of bounded semi-algebraically connected components of the realizable
sign conditions of ${\mathcal P}$. 
It follows from Lemma \ref{lem:2} that if $C$ is a  
bounded semi-algebraically connected component of $\RR(\sigma,\R^k)$ for 
$\sigma \in {\rm Sign}({\mathcal P}),$ 
then $b^{BM}_{\ell}(C) \geq 1.$

It now follows from Lemma \ref{lem:3} and Inequality \eqref{eqn:proofmain5}
that the number of bounded semi-algebraically connected 
components of $\RR(\sigma,\R^k)$ over all 
$\sigma \in {\rm Sign}({\mathcal P})$,  
is bounded by
\begin{eqnarray*}
\sum_{{\mathcal Q} \subset {\mathcal P}} b^{BM}(W_{{\mathcal Q},\R})
&\leq&
\sum_{0 \leq \ell \leq k}\left({s \choose k - \ell}{s  \choose \ell}d^k + 
\sum_{1 \leq i \leq \ell} {s \choose k - \ell}{s \choose \ell - i}d^{O(k^2)}
\right)\\
&\leq &
\sum_{0 \leq \ell \leq k} {s \choose k - \ell}{s  \choose \ell}d^k + 
O(s^{k-1}) \\
&=&
\left(\sum_{\ell=0}^k \frac{d^k}{\ell! (k -\ell)!}\right)s^k + O(s^{k-1}) \\
&=& \frac{(2d)^k}{k!} s^k + O(s^{k-1}).
\end{eqnarray*}

In order to bound the number of unbounded components, consider the 
realizations of sign conditions of the family ${\mathcal P}$ restricted
to the $(k-1)$-dimensional real variety $\partial \overline{B_{k}(0,\Omega)}$. 
The number of unbounded semi-algebraically connected components of the realizable
sign conditions of ${\mathcal P}$ over $\R^k$, is bounded by the
number of semi-algebraically connected components of the realizations of sign conditions
of ${\mathcal P}$ restricted to the variety $\partial \overline{B_{k}(0,\Omega)}$.
It follows from the bound in (\ref{eqn:pams}) that
the number of semi-algebraically connected components of the realizations of sign conditions
of ${\mathcal P}$ restricted to the variety $\partial \overline{B_{k}(0,\Omega)}$
is bounded by $O(s^{k-1})$. Thus,
the number of unbounded semi-algebraically connected components 
of the realizable sign conditions of ${\mathcal P}$ over $\R^k$
is also bounded by $O(s^{k-1})$. 
This proves the theorem.
\end{proof}

\begin{proof}[Proof of Theorem \ref{the:0-1}]
Let $\rho \in \mbox{Zero-pattern}({\mathcal P})$ with $d_\rho = \ell.$
Then there exists an irreducible component, 
$V_{\rho,i}$ of $V_\rho$, with $\dim_{\C} V_{\rho,i} = \ell$,
and such that $\dim_{\C} (V_{\rho,i} \cap \RR(\rho,\C^k))$ is also equal to
$\ell$. Moreover, for any $\rho' \in \mbox{Zero-pattern}({\mathcal P}),$
with $\rho' \neq \rho,$ we must have that
$\dim_{\C} (V_{\rho,i} \cap \RR(\rho',\C^k)) < \ell,$
since $\RR(\rho,\C^k) \cap \RR(\rho',\C^k) = \emptyset$ and
$\dim_{\C}(V_{\rho,i} \setminus \RR(\rho,\C^k))$ is clearly
$< \ell$.

Thus, if we charge $\rho$ to $V_{\rho,i}$,
it is clear that $V_{\rho,i}$
cannot be charged by any zero pattern other than $\rho$. 
It now follows from Lemma \ref{lem:positive}
that the number of distinct $V_{\rho,i}$ of dimension $\ell$ is bounded by 
${s \choose k-\ell} d^{k - \ell}.$
\end{proof}

\bibliographystyle{amsplain}
\bibliography{master}

\end{document}